\documentclass[a4paper,10pt]{article}
\usepackage{amssymb,graphicx,graphics,amsmath}
\usepackage{makeidx}
\makeindex
\newtheorem{next}{}[section]
\newcommand{\nx}[1]{{\bf #1.}\ }
\newcommand{\proof}{\indent{\it Proof.}\quad}
\newcommand{\qed}{\hfill\rule{1.6mm}{1.6mm}\,}

\renewcommand {\k} {\Bbbk}

\newcommand {\Z} {{\mathbb Z}}
\newcommand {\N} {{\mathbb N}}
\newcommand {\Q} {{\mathbb Q}}

\newcommand{\si}{\k[{\bf x}]}
\newcommand{\sis}{\k[{\bf z}]}

\newcommand{\M}{{\cal M}}

\newcommand{\I}{{\cal I}}

\newcommand{\dtb}{$\stackrel{\mbox{\circle*{4.0}}}{}$}

\newcommand{\igno} [1] {}

\title{Sums of toric ideals}
\author{DE ALBA CASILLAS Hernan, MORALES Marcel}
\date{}

\begin{document}
\maketitle

{\small \abstract{Given two toric ideals $I_1,I_2\subset\si$, it is not always true that $I_1+I_2$ is a toric ideal. Given $I_1,\dots,I_k\subset\si$ a familly of toric ideals we give necessary conditions in order to have that $I_1+\dots+I_k$ is a toric ideal.}
\section {Introduction}
Let $\k$ be a field. Let $\k[\bf{t}^{\pm}]:=\k[t_1^{\pm},\dots,t_m^{\pm}]$ the Laurent polynomial ring on $m$ variables over a field $\k$. Let $\alpha=(\alpha^{(1)},\dots,\alpha^{(m)})\in \Z^{m}$, ${\bf t}^{\alpha}=t_1^{\alpha^{(1)}}\dots t_{m}^{\alpha^{(m)}}$ is called a monomial of $\k[\bf{t}]$. Let $Y_1={\bf y}^{\alpha_1},\dots,Y_n={\bf y}^{\alpha_n}$ be monomials in $\k[\bf{t}]$, we denote by $S:=\k[Y_1,\dots Y_m]$  the subalgebra of the polynomial ring over the field $k$ generated by $Y_1,\dots,Y_m$. The set of monomials $\M$ contained in $S$ form a semigroup under multiplication, and the function ${\rm deg}:\M\rightarrow \N^m$ that assigns each monomial its exponent vector maps $\M$ isomorphically onto a subsemigroup $H$ of $\N^m$. Up to isomorphism, $S$ is the semigroup algebra $\k[H]$. We will always assume that $Y_1,\dots,Y_n$ are irreducible elements of $\M$, or, in other words, that they form a minimal system of algebra generators of $S$.\\

Now, let $\si:=\k[x_1,\dots,x_n]$ the polynomial ring on $n$-variables and $\varphi:\si\rightarrow S$ the algebra morphism defined by $\varphi(x_i)=Y_i$. As $S$ is a Laurent polynomial ring, $S$ is an integer domain and $J={\rm Ker}\mbox{ }\varphi$ is prime. We will say that $J$ is a {\bf toric ideal} of $\si$ and $\varphi$ is a parametrization of $J$. The next two propositions are well known:

\begin{next}\nx{Proposition}\label{remarque-def-torique}
 Let $A=[\alpha_1,\dots,\alpha_n]\in M_{m,n}(\Z)$, $\pi:\N^n\rightarrow \Z^m$ the semigroup morphims defined by $\pi({\bf e}_i)=\alpha_i$ and $\varphi:\si\rightarrow \k[t^{\alpha_1},\dots,t^{\alpha_n}]$ the algebra morphism defined by $\varphi(x_i)=t^{\alpha_i}$. Then ${\rm ker}\mbox{ }\varphi=I_A:=({\bf x}^{u_+}-{\bf x}^{u_{-}}:\pi(u_+)=\pi(u_{-})).$ Moreover, for any $u\in\N^n$, $\varphi(z^u)=t^{A\cdot u}$.
\end{next}

\begin{next}\nx{Proposition}\label{dim-rang}
 Let $\varphi:\si\rightarrow\k[{\bf t}^{\pm 1}]$ be a parametrization of $I:={\rm Ker}\mbox{ }\varphi$ such that $\varphi(x_i)=t^{\alpha_i}$. Let set $A=[\alpha_1,\dots,\alpha_n]\in M_{m,n}(\Z)$, then ${\rm rank}\mbox{ }(A)={\rm dim}\mbox{ }(\si/I)$.
\end{next}
\igno{\proof
D'apr\`es le point $6$ de la proposition \ref{torique-equiv} $\si/I\cong\si[{\bf t}^{\alpha_1},\dots,{\bf t}^{\alpha_n}]$. La dimension de Krull du domaine entier $\k[{\bf t}^{\alpha_1},\dots,{\bf t}^{\alpha_n}]$ est le plus grand nombre de mon\^omes ${\bf t}^\alpha_i$ alg\'ebriquement ind\'ependants. Or, un ensemble de mon\^omes est alg\'ebriquement ind\'ependant si et seulement si ses vecteurs exposants sont lin\'eairement ind\'ependants. Ainsi ${\rm rang}\mbox{ }(A)={\rm dim}\mbox{ }(\si/I)$.
\qed}

\begin{next}\nx{Definition}\label{rang-maximal}
 {\rm Let $\varphi:\si\rightarrow\k[{\bf t}^{\pm 1}]$ a parametrization of $I:={\rm Ker}\mbox{ }\varphi$. We say that $\varphi$ has maximal rank if ${\rm dim}\mbox{ }(\si/I)$ is equal to cardinality of ${\bf t}$.}
\end{next}

We say that $u=(u_1,\dots,u_n)\in\N^n$ and $v=(v_1,\dots,v_n)\in\N^n$ have disjoint support if for all $i\in\N$ such that $u_i\neq 0$ then $v_i=0$ and reciprocally. Furthermore, for any $u\in\Z^n$ there is two unique vectors $u_+,u_{-}\in\N^n$ with disjoint supports such that $u=u_+-u_{-}$. In this manner, by proposition \ref{remarque-def-torique} we deduce the next proposition:

\begin{next}\nx{Proposition}\label{pour-changement-de-base}
Let $\varphi:\si\rightarrow\k[{\bf t}^{\pm 1}]$ be a parametrization of ${\rm Ker}\mbox{ }\varphi$. Then, setting $${\rm Ker}_{\Z}\mbox{ }A:=\{u\in\Z^n:Au=0\},$$ ${\rm Ker}\mbox{ }\varphi=({\bf x}^{u_+}-{\bf x}^{u_{-}}:u_+,u_{-}\in\N^n\mbox{ }{\rm of}\mbox{ }{\rm disjoint}\mbox{ }{\rm support}\mbox{ }{\rm and}\mbox{ }(u_{+}-u_{-})\in{\rm Ker}_{\Z}\mbox{ }A)$.
\qed
\end{next}

From the proposition \ref{pour-changement-de-base}, in order to describe a toric ideal $I$, it is enough to solve the linear homogene equation system over $\Z$ $$A\cdot\left(\begin{array}{l}u_1\\ \vdots\\ u_n \end{array}\right)=0,$$ where $A$ is a matrix related to the parametrization $\varphi$ of $I$. We know from linear algebra that for any $Q\in M_{m,m}(\Q)$ non-singular, ${\rm Ker}_{\Z}\mbox{ }A={\rm Ker}_{\Z}\mbox{ }QA$, thus the next proposition is immediate.

\begin{next}\nx{Proposition}\label{changement de base}
 If $\widetilde{A}\in M_{m,n}(\Z)$ is a matrix obtained from $A$ by a left multiplication for a non-singular matrix $Q\in M_{m,m}(\Q)$, then $I_{\widetilde{A}}={\rm Ker}\mbox{ }\varphi$ and the matrix $\widetilde{A}$ define a new parametrization of ${\rm Ker}\mbox{ }\varphi$ par $\widetilde{\varphi}(x^u)={\bf t}^{[\widetilde{A}]\cdot u}$. In particular, if $\widetilde{A}\in M_{m,n}(\Z)$ is obtained by a finite sequence of elementary operations over the lines of $A$, then $I_{\widetilde{A}}={\rm Ker}\mbox{ }\varphi$ and the matrix $B$ define a new parametrization of ${\rm Ker}\mbox{ }\varphi$. Reciprocally, if $I\subset\k[x_1,\dots,x_n]$ is a toric ideal, then $L=\{u\in \Z^n:x^{u_+}-x^{u_-}\in I\}$ is a saturated lattice. If $\{b_1,\dots,b_m\}$ is a basis of $L$ and $B=[b_1\dots b_m]\in M_{n,m}$ then any matrix $A\in M_{m,n}(\Z)$ such that $AB=0$, $A$ is a parametrization of $I$.
\end{next}
\igno{\proof
Due to \ref{pour-changement-de-base} ${\rm Ker}_{\Z}\mbox{ }A={\rm Ker}_{\Z}\mbox{ }B$. Thus ${\rm Ker}\mbox{ }\varphi={\rm Ker}\mbox{ }\widetilde{\varphi}$.}
\qed
\section {Sums of toric ideals}

\begin{next}\nx{D\'efinition}\label{homog-des}
 {\rm We consider two polynomial rings $\k[{\bf z}]:=\k[z_1,\dots,z_n]$ and $\k[x][z]\cong\k[x,z_1,\dots,z_n]$. Let $f:=f(z_1,\dots,z_n)$ be a polynomial in $\sis$ and $g:=g(z_1,\dots z_n,x)$ a homogeneous polynomial in $\si[{\bf z}]$. Writting $f=\sum_{i=1}^m a_i{\bf z}^\beta_i$ for some $a_i\in\k$ and $\beta_i\in \N^n$, we define the degree of $f$ as\break ${\rm deg}_{\N}(f)={\rm max}_{i=1}^{n}(|\beta_i|)$. The {\bf homogenization} of $f$ by $x$ is the polynomial $f^{{\rm hom}_{x}}=x^{{\rm deg}_{\N}(f)}f(\frac{z_1}{x},\dots,\frac{z_n}{x})$. The {\bf dehomogenization} of $g$ by $x$ is the polynomial $g^{{\rm deh}_x}=g(z_1,\dots z_n,1)\in\sis$. In case where there is not confusion about the variable $x$, we will only write $f^{{\rm hom}}$ and $g^{{\rm deh}}$. In addition, let $I$ be an ideal of $\sis$ and $J$ an ideal of $\k[{\bf z},x]$. The homogenization of $I$ by $x$ is the ideal $I^{{\rm hom}}=(f^{\rm hom}:f\in I)$ and the dehomogenization of $J$ by $x$ is the ideal $J^{{\rm deh}}=(g^{\rm deh}:g\in J)$.}
\end{next}

\begin{next}\nx{Lemma}\label{lemme-homog}{\rm \cite[\rm Lemma 4.14]{sturmfels}}
 Let $A=[\alpha_1,\dots,\alpha_n]\in M_{m,n}(\Z)$. The ideal $I_A$ is homogeneous if and only if there exists a vector $\omega\in\Q^m$ such that $\alpha_i\cdot\omega=1$ for any $i=1,\dots,n$ for which $\alpha_i\neq 0$.
\end{next}

\begin{next}\nx{Lemma}\label{torique-homogene}
 Let ${\bf z}=(z_1,\dots,z_n)$ be a sequence of variables, $x$ be a variable which does not belong to ${\bf z}$ and $I\subset \k[{\bf z},x]$ be a homogeneous toric ideal such that there exist a parametrization of $I$:
\begin{center}
$\varphi:\k[{\bf z},x]\rightarrow\k[{\bf t}^{\pm 1},s]$
\end{center}
with $\varphi(z_i)={\bf t}^{\alpha_i}{s}^{\gamma_i}$ where $\alpha_i\in \Z^{n}$, $\gamma_i\in\Z$, and $\varphi(x)=s^{\gamma_{n+1}}$ where $\gamma_{n+1}\in\Z^*$. Let $\widetilde{\varphi}:\sis\rightarrow \k[{\bf t}^{\pm 1}]$ be the morphism defined by $\widetilde{\varphi}(z_i)=(\varphi(z_i))^{{\rm deh}_s}$, then $({\rm Ker}\mbox{ }\widetilde{\varphi})^{\rm hom}=I$.
\end{next}
\proof
We set $\widetilde{A}=[\alpha_1,\dots,\alpha_n]\in M_{m,n}(\Z)$, $\gamma=(\gamma_1,\dots,\gamma_n)$ and
$$A=\left(\begin{array}{c|c}
\widetilde{A}&{\bf 0}\\ \hline
\gamma&\gamma_{n+1}
\end{array}\right)\in M_{m,n}(\Z) \mbox{ }{\rm so}\mbox{ }{\rm that}\mbox{ } I:={\rm Ker}\mbox{ }\varphi=I_{A}.$$ 
Now, we are going to prove that  $({\rm Ker}\mbox{ }\widetilde{\varphi})^{\rm hom}=I$. Clearly $I\subset ({\rm Ker}\mbox{ }\widetilde{\varphi})^{\rm hom}$; so we only need to prove that  $({\rm Ker}\mbox{ }\widetilde{\varphi})^{\rm hom}\subset I$.\\

Let ${\bf z}^{u_+}-{\bf z}^{u_{-}}\in {\rm Ker}\mbox{ }\widetilde{\varphi}$, where $u_{+}=(u_{+}^{(1)},\dots,u_{+}^{(n)})$ and $u_{-}=(u_{-}^{(1)},\dots,u_{-}^{(n)})\in \N^n$ have disjoint support. We should prove that $({\bf z}^{u_+}-{\bf z}^{u_{-}})^{\rm hom}\in I$. We recall that if $u=(u^{(1)},\dots,u^{(n)})\in\Z^n$, $|u|=\sum_{i=1}^n u^{(i)}$. Without lost of generality we can assume $|u_+|\leq|u_{-}|$. Setting $u_{+}^{(n+1)}\in\N$ such that $|u_+|+u_{+}^{(n+1)}=|u_{-}|$, we have that $$({\bf z}^{u_+}-{\bf z}^{u_{-}})^{\rm hom}={\bf z}^{u_+}x^{u_{+}^{(n+1)}}-{\bf z}^{u_{-}}.$$
Let set $\widetilde{u}_{+}=(u_{+},u_{+}^{(n+1)})\in\N^{n+1}$, and $\widetilde{u}_{-}=(u_{-},0)\in\N^{n+1}$. Since $I$ is homogeneous, thanks to lemma \ref{lemme-homog} there exists $\omega\in\Q^m$ and $\omega_{m+1}\in \Q$ such that $(\omega,\omega_{m+1})\cdot\alpha'_i=1$ for any $i=1,\dots,n+1$, this is equivalent to $$(\omega,\omega_{m+1})\cdot A=(1,\dots,1).$$
Thus, $$\begin{array}{ll}
	|u_{+}|+u_{+}^{(n+1)}&=(1,\cdots 1)\cdot\widetilde{u}_{+}\\
	&=((\omega,\omega_{m+1})\cdot A )\cdot(u_{+},u_{+}^{(n+1)})\\
	&=\omega\cdot (\widetilde{A}\cdot u_{+})+\omega_{m+1}((\gamma,\gamma_{n+1})\cdot \tilde{u}_{+})\\
	&=\omega\cdot (\widetilde{A}\cdot u_{+})+\omega_{m+1}\sum_{i=1}^{n+1}\gamma_i u_{+}^{(i)},
        \end{array}$$
and
$$\begin{array}{ll}
|u_{-}|&=(1,\cdots 1)\cdot\widetilde{u}_{-}\\
&=((\omega,\omega_{m+1})\cdot A) \cdot(u_{-},0)\\
&=\omega\cdot (\widetilde{A}\cdot u_{-})+\omega_{m+1}\sum_{i=1}^{n}\gamma_i u_{-}^{(i)}.
\end{array}$$
Besides $0=\widetilde{\varphi}(\bf{z}^{u_+}-{\bf z}^{u_{-}})={\bf t}^{\widetilde{A}\cdot u_{+}}-{\bf t}^{\widetilde{A}\cdot u_{-}}$, so $\widetilde{A}\cdot u_{+}=\widetilde{A}\cdot u_{-}$. Furthermore, $|u_+|+u_{n+1}=|u_{-}|$, and then $\omega_{m+1}\sum_{i=1}^{n+1}\gamma_i u_{+}^{(i)}=\omega_{m+1}\sum_{i=1}^{n}\gamma_i u_{-}^{(i)}$. In addition, $\omega_{m+1}\neq 0$, because $$1=(\omega,\omega_{m+1})\cdot\alpha'_{m+1}=\omega_{m+1}\gamma_{m+1}.$$ Thus, 
\begin{equation}\label{a-voir-zero}
\sum_{i=1}^{n+1}\gamma_i u_{+}^{(i)}=\sum_{i=1}^{n}\gamma_i u_{-}^{(i)}.
\end{equation}
On the other side
$$\begin{array}{ll}
\varphi(\bf{z}^{u_+}x^{u_{+}^{(n+1)}}-{\bf z}^{u_{-}})&=\varphi(z^{u_+}x^{u_{+}^{(n+1)}})-\varphi({\bf z}^{u_{-}})\\
&=({\bf t},s)^{A\cdot u_{+}'}-({\bf t},s)^{A\cdot u_{-}'}\\
&={\bf t}^{\widetilde{A}\cdot u_{+}}s^{\gamma\cdot u_{+}} s^{\gamma_{n+1} u_{+}^{(n+1)}}-{\bf t}^{\widetilde{A}\cdot u_{-}}s^{\gamma\cdot u_{-}}\\
&={\bf t}^{\widetilde{A}\cdot u_{+}}s^{\sum_{i=1}^{n+1}\gamma_i u_{+}^{(i)}}-{\bf t}^{\widetilde{A}\cdot u_{-}}s^{\sum_{i=1}^{n}\gamma_i u_{-}^{(i)}}\\
&={\bf t}^{\widetilde{A}\cdot u_{+}}s^{\sum_{i=1}^{n+1}\gamma_i u_{+}^{(i)}}-{\bf t}^{\widetilde{A}\cdot u_{+}}s^{\sum_{i=1}^{n}\gamma_i u_{-}^{(i)}}\\
&={\bf t}^{\widetilde{A}\cdot u_{+}}(s^{\sum_{i=1}^{n+1}\gamma_i u_{+}^{(i)}}-s^{\sum_{i=1}^{n}\gamma_i u_{-}^{(i)}})=0.
\end{array}$$
The last inequality is due to the equation \ref{a-voir-zero}. By consequence
$$({\bf z}^{u_+}-{\bf z}^{u_{-}})^{\rm hom}={\bf z}^{u_+}x^{u_{+}^{(n+1)}}-{\bf z}^{u_{-}}\in I$$ and $({\rm Ker}\mbox{ }\widetilde{\varphi})^{\rm hom}\subset I$. As we have already seen $I\subset({\rm Ker}\mbox{ }\widetilde{\varphi})^{\rm hom}$, we conclude $({\rm Ker}\mbox{ }\widetilde{\varphi})^{\rm hom}=I.$

\qed

\begin{next}\nx{Lemma}\label{inversible}
Let $A\in M_{m,n}(\Z)$ where $m\leq n$ and ${\rm rang}\mbox{ }(A)=m$. Then, for any $i\in\{1,\dots,n\}$ such that the $i$-th column of $A$ is not nul, there always exists a submatrix $A'$ of $m$ columns of $A$ where the $i$-th column of $A$ is a column of $A'$ and $A'$ is non-singular over $\Q$.
\end{next}
\proof
Since ${\rm rank}\mbox{ }(A)=m$, there exists a matrix $A_1\in M_{m,m}(\Z)$ with $m$ columns of $A$ such that $A_1$ is inversible. Let $\alpha_i$ be the $i$-th column of $A$ with $\alpha_i\neq 0$ and $\alpha_1',\dots,\alpha_m'\in \Z^m$ the $m$-columns of $A_1$. Then $\{\alpha_1',\dots,\alpha_m'\}$ is a basis of $\Q^m$. Thus, there exists $a_1,\dots,a_m\in \Q$ not all nul such that $\alpha_i=\sum_{j=1}^m a_j\alpha_j'$. We can assume $\alpha_1\neq 0$, then $\alpha_1'=\sum_{j=2}^m \frac{a_j}{\alpha_1}\alpha_j'-\frac{1}{a_1}a_i$. So $\{\alpha_i,\alpha_2',\dots,\alpha'_m\}$ is a basis of $\Q^m$ and the matrix $A'=[\alpha_i,\alpha_2',\dots,\alpha_m']$ is non-singular over $\Q$.

\qed 

\begin{next}\nx{Lemma}\label{lemme-toriquin}
Let $\varphi:\si\rightarrow\k[t_1,\dots,t_{m'},t_1^{-1},\dots,t_{m'}^{-1}]$ be a parametrization of the ideal $I:={\rm Ker}\mbox{ }\varphi$ and $A^*\in M_{m',n}(\Z)$ the matrix related to this parametrization such that $${\rm rang}\mbox{ }(A^*)=m\leq m'.$$ For any $i\in\{1,\dots,n\}$ such that $\varphi(x_i)\neq 1$, there exists a parametrization of $I$ $$\widetilde{\varphi}:\si\rightarrow\k[t_1,\dots,t_{m},t_1^{-1},\dots,t_{m}^{-1}]$$ of maximal rank such that $\widetilde{\varphi}(x_i)=t_j^q$, for some $q\in\N$ and $j\in\{1,\dots m\}$ .
\end{next}
\proof
By hypothesis $\varphi(x_i)=t_1^{\alpha_{i,1}}\cdots t_{m'}^{\alpha_{i,m'}}\neq 1$. Let $j'\in\{1,\dots,m'\}$ be the smallest integer such that $\alpha_{i,j'}\neq 0$. Let denote by $\beta_r$ the $r$-th line vector  of the matrix $A^*$. As ${\rm rang}(A^*)=m$, there exist $\beta_{i_1},\dots,\beta_{i_m}\in \{\beta_1,\dots,\beta_{m'}\}$ such that $${\rm Span}_{\Q}(\beta_{i_1},\dots,\beta_{i_m})\cong \Q^m.$$ Thus there exist $b_1,\dots,b_m\in \Q$ not all null such that $\beta_{j'}=\sum_{k=1}^m b_{k}\beta_{i_k}$. We can assume $b_1\neq 0$, then $$\beta_{i_1}=\sum_{k=2}^m \frac{b_k}{b_1}\beta_{i_k}-\frac{1}{b_1}\beta_{j'}.$$ So $\beta_{j'},\beta_{i2},\dots,\beta_{im}$ is a basis of ${\rm Span}_{\Q}(\beta_{i_1},\dots,\beta_{i_m})=\Q^m$ and there exists a non-singular matrix $B\in M_{m',m'}(Q)$ such that 

\begin{center}
$B\cdot A^*=\left(\begin{array}{c}A\\{\bf 0}\end{array}\right)\in M_{m',n}(\Q)$, avec $A=(\beta_{j'},\beta_{i_{k_2}},\dots,\beta_{i_{k_m}})^{T}$ et ${\bf 0}\in M_{m'-m,n}(\Q)$.
\end{center} 

By the lemma \ref{inversible} there exists a non-singular matrix $A'\in M_{m,m}(\Z)$ of m columns of $A$ containning the $i$-th column of $A$. Moreover, we can assume that
$A=\left[
 A'|A''\right]$
where $A''\in M_{m,n-m}$. Since $A'$ is invertible, there exists $B\in M_{m,m}(\Q)$ such that $BA'=I$. By this way $BA=[BA'|BA'']=[I|BA'']$. Let $q$ be the smaller natural integer such that $$q(BA)=(qB)A=[qI|q(BA'')]\in M_{m,n}(\Z).$$ From the  proposition \ref{changement de base}, $(qB)A$ defines a parametrization of $I$, $\widetilde{\varphi}:\si\rightarrow\k[{\bf t}^{\pm 1}]$ of maximal rank (${\rm rank}\mbox{ }(A')=m$)  such that there exists $j\in\{1,\dots m\}$ and $\widetilde{\varphi}(x_i)=t_j^q$.
\qed

\begin{next}\nx{Theorem}\label{sum-2-torique}
 Let ${\bf z_{1}}=\{z_{1,1},\dots,z_{1,n_1}\}$ and ${\bf z_{2}}=\{z_{2,1},\dots,z_{2,n_2}\}$ be two disjoint variable sets, $x$ be a new variable and $I_1\subset \k[{\bf z_{1}},x]$, $I_2\subset\k[{\bf z_{2}},x]$ two homogeneous toric ideals such that $x$ appears in at least one generator of a minimal generator system of $I_1$ and $I_2$. Then:

\begin{enumerate}
\item $I_1+I_2$ is a homogeneous toric ideal of $\k[{\bf z_{1}},{\bf z_{2}},x]$.
\item ${\rm dim}\mbox{ }(\k[{\bf z_{1}},{\bf z_{2}},x]/(I_1+I_2))={\rm dim}\mbox{ } (\k[{\bf z_{1}},x]/I_1)+{\rm dim}\mbox{ }(\k[{\bf z_{2}},x]/I_2)-1$.
\end{enumerate}
\end{next}
\proof
Without lost of generality we can suppose $I_1$ and $I_2$ are non degenerated, i.e., for any $i=1,2$, any variable of ${\bf z_i}\cup \{x\}$ appears in a minimal system of generators of $I_i$.\\

Before to prove this theorem, we make the following remarks:\\
We set ${\bf z}={\bf z_{1}}\cup {\bf z_{2}}$. From the lemma \ref{lemme-toriquin} there exist
$$\overline{\varphi}_1:\k[{\bf z_1},x]\rightarrow\k[{\bf t}^{\pm 1},s^{\pm 1}]\mbox{ }{\rm et}\mbox{ }\overline{\varphi_2}:\k[{\bf z_1},x]\rightarrow\k[{\bf w}^{\pm 1},s^{\pm 1}]$$ 
parametrizations of maximal rank of $I_1$ and $I_2$ respectively, where $\overline{\varphi}_i(x) =s^{\gamma_{i,n_i+1}}$ and\break $\gamma_{i,n_i+1}\in \Z^*$ for any $i=1,2$. We set $m_1$, $m_2$ the cardinals of ${\bf t}$ and ${\bf u}$, repectively.  We set also $\gamma={\rm lcm}( \gamma_{1,n_1+1}, \gamma_{2,n_2+1})$. We denote by $\overline{A}_1\in M_{m_1+1,n_1+1}(\Z)$ et $\overline{A}_2\in M_{m_2+1,n_2+1}(\Z)$ the matrices which represent the parametrizations $\overline{\varphi_1}$ and $\overline{\varphi_2}$, respectively. These matrices can be written as below:
$$\overline{A}_1=\left(\begin{array}{c|c}
\overline{A}_1'&{\bf 0}\\ \hline
\alpha_1'&\gamma_{1,n_1+1}
\end{array}\right)\mbox{ }{\rm et} 
\overline{A}_2=\left(\begin{array}{c|c}
\overline{A}_2'&{\bf 0}\\ \hline
\alpha_2'&\gamma_{2,n_2+1}
\end{array}\right),$$
where $\overline{A}_1'\in M_{m_1,n_1}(\Z)$, $\overline{A}_2\in M_{m_2,n_2}(\Z)$, $\alpha_1=(\gamma_{1,1},\dots,\gamma_{1,n_1})\in\Z^{n_1}$ and $\alpha_2=(\gamma_{2,1},\dots,\gamma_{2,n_2})\in\Z^{n_2}$. Multiplying the last line of each matrix $\overline{A}_1$ and $\overline{A}_2$ respectively, by $\frac{\gamma}{\gamma_{1,n_1+1}}$ and $\frac{\gamma}{\gamma_{2,n_2+1}}$, where $\gamma=lcm(\gamma_{1,n_1+1},\gamma_{2,n_2+1})$, we obtain from the proposition \ref{changement de base} other parametrizations $\varphi_1$ and $\varphi_2$ for $I_1$ and $I_2$, respectively, namely
$$A_1=\left(\begin{array}{c|c}
A_1'&{\bf 0}\\ \hline
\alpha_1&\gamma
\end{array}\right)\mbox{ }{\rm et} 
A_2=\left(\begin{array}{c|c}
A_2'&{\bf 0}\\ \hline
\alpha_2&\gamma
\end{array}\right),$$
repectively, where $\alpha_1=\frac{\gamma}{\gamma_{1,n_1+1}}\cdot\alpha_1'\in\Z^{n_1}$ et $\alpha_1=\frac{\gamma}{\gamma_{1,n_2+1}}\cdot\alpha_2'\in\Z^{n_2}$. Let $\overline{\varphi}:\k[{\bf z},x]\rightarrow\k[{\bf t}^{\pm 1},{\bf u}^{\pm 1},s^{\pm 1}]$ be the morphism defined by $\overline{\varphi}(x)=s^\gamma$ and $\overline{\varphi}(z_{j,i})=\varphi_i(z_{j,i})$ o\`u $i=1,2$.\\

Now we are able to prove the theorem:
 
\begin{enumerate}
 \item We set $J:={\rm Ker}\mbox{ }\overline{\varphi}$, then $J$ is a toric ideal and we have that:
\begin{itemize}
 \item [\dtb] $J$ is homogeneous. This is obtained by the following reasoning:\\
Since  $I_1$ and $I_2$ are homegeneous ideal, by the lemma \ref{lemme-homog}, there exist $\omega_1\in\Q^{m_1}$, $\omega_2\in\Q^{m_2}$, $\omega'\in\Q$ and $\omega'\in\Q$ such that

$$(1,\dots,1)=(\omega_1,\omega')\cdot A_1=(\omega_1\cdot A_1'+\omega'\cdot\alpha_1,\omega'\cdot\gamma)\mbox{ }  {\rm et}$$ $$(1,\dots,1)=(\omega_2,\omega'')\cdot A_2=(\omega_2\cdot A_2'+\omega''\cdot\alpha_2,\omega''\cdot\gamma).$$
Thus $\omega'\cdot\gamma=\omega''\cdot\gamma$ and by consequence $\omega'=\omega''$
Moreover, the matrix which represents the morphism $\overline{\varphi}$ is
$$A=\left(\begin{array}{c|c|c}
A_1'&{\bf 0}&{\bf 0}\\ \hline
{\bf 0}&A_2'&{\bf 0}\\ \hline
\alpha_1&\alpha_2&\gamma
\end{array}\right)$$ and
$$\begin{array}{ll}(\omega_1,\omega_2,\omega')\cdot A&=(\omega_1\cdot A_1+\omega'\cdot\alpha_1,\omega_2\cdot A_2+\omega'\cdot\alpha_2,\omega'\gamma)\\
&=(\omega_1\cdot A_1+\omega'\cdot\alpha_1,\omega_2\cdot A_2+\omega'\cdot\alpha_2,\omega'\gamma)\\
&=(1,\dots,1,1,\dots,1,1).
  \end{array}$$
Due to lemma \ref{lemme-homog}, $J={\rm Ker}\mbox{ }\overline{\varphi}=I_A$ is homogeneous.

 \item [\dtb] Let us prove that $J=I_1+I_2$. It is clear that $I_1+I_2\subset J$. We need to prove that $J\subset I_1+I_2$:\\
Let $B:={\bf z_1}^{u_1}{\bf z_2}^{u_2}x^\delta-{\bf z_1}^{v_{1}}{\bf z_2}^{v_{2}}\in J$ homogeneous, where $\delta\in\N$ and $u_{i},v_{i}\in \N^{n_i}$ for any $i=1,2$, with the property that if $(u_{i})_j\neq 0$, then $(v_{i})_j=0$ and reciprocally. Thus
\begin{equation}\label{deg-hom}
 |u_1|+|u_2|+\delta={\rm deg}_{\N}({\bf z_1}^{u_1}{\bf z_2}^{u_2}x^\delta)={\rm deg}_{\N}({\bf z_1}^{v_{1}}{\bf z_2}^{v_{2}})=|v_1|+|v_2|.
\end{equation}
Furthermore, since $$0=\overline{\varphi}(B)\widetilde{\varphi_1}({\bf z_1}^{u_1})\widetilde{\varphi_2}({\bf z_2}^{u_2})s^{\lambda_1}-\widetilde{\varphi_1}({\bf z_1}^{v_1})\widetilde{\varphi_2}({\bf z_2}^{v_2})s^{\lambda_2}$$ and the variable sets ${\bf t}$, ${\bf w}$ and $s$ are disjoint,
\begin{center}
 $\widetilde{\varphi_1}({\bf z_1}^{u_1})=\widetilde{\varphi_1}({\bf z_1}^{v_1})$ et $\widetilde{\varphi_2}({\bf z_2}^{u_2})=\widetilde{\varphi_2}({\bf z_2}^{v_2})$ et $\lambda_1=\lambda_2$.
\end{center}
So ${\bf z_1}^{u_1}-{\bf z_1}^{v_1}\in({\rm Ker}\mbox{ }\widetilde{\varphi_1})$ et ${\bf z_2}^{u_2}-{\bf z_2}^{v_2}\in({\rm Ker}\mbox{ }\widetilde{\varphi_2})$. By the lemma \ref{torique-homogene}: 
\begin{center}
 $({\bf z_1}^{u_1}-{\bf z_1}^{v_1})^{\rm hom}\in({\rm Ker}\mbox{ }\widetilde{\varphi_1})^{\rm hom}=I_1$ et $({\bf z_2}^{u_2}-{\bf z_2}^{v_2})^{\rm hom}\in({\rm Ker}\mbox{ }\widetilde{\varphi_2})^{\rm hom}=I_2$.
\end{center}

\begin{itemize}
 \item [\dtb] We will prove that $B$ is a algebric combination of $({\bf z_1}^{u_1}-{\bf z_1}^{v_1})^{\rm hom}$ and $({\bf z_2}^{u_2}-{\bf z_2}^{v_2})^{\rm hom}$. In order to prove it, we consider three cases:
	\begin{enumerate}
	 \item Let $|u_1|<|v_1|$. Then, there exists $\alpha\in\N$, $\alpha>0$, such that $$({\bf z_1}^{u_1}-{\bf z_1}^{v_1})^{\rm hom}={\bf z_1}^{u_1}x^{\alpha}-{\bf z_1}^{v_1}.$$
So 
	\begin{equation}\label{eq-aux1}
	|u_1|+\alpha=|v_1|.
	\end{equation}
Now, we consider the following two subcases:
		\begin{enumerate}
			\item Let $|u_2|\leq |v_2|$. Then, there exist $\beta\in\N$, such that $$({\bf z_2}^{u_2}-{\bf z_2}^{v_2})^{\rm hom}={\bf z_2}^{u_2}x^\beta-{\bf z_2}^{v_2}.$$ Then $|u_2|+\beta=|v_2|$. Furthermore, from the equations \ref{eq-aux1} and \ref{deg-hom}, we deduce that 
			$$(|u_1|+\alpha)+(|u_2|+\beta)=|v_1|+|v_2|=|u_1|+|u_2|+\delta.$$
			By consequence $\alpha+\beta=\delta$ and
			$$\begin{array}{ll}
			B&={\bf z_1}^{u_1}{\bf z_2}^{u_2}x^\delta-{\bf z_1}^{v_{1}}{\bf z_2}^{v_{2}}\\
			&={\bf z_1}^{u_1}x^\alpha({\bf z_2}^{u_2}x^\beta-{\bf z_2}^{v_2})+{\bf z_2}^{v_2}({\bf z_1}^{u_1}x^{\alpha}-{\bf z_1}^{v_1})\\
			&={\bf z_1}^{u_1}x^\alpha({\bf z_2}^{u_2}-{\bf z_2}^{v_2})^{\rm hom}+{\bf z_2}^{v_2}({\bf z_1}^{u_1}-{\bf z_1}^{v_1})^{{\rm hom}}\in I_2+I_1.
			\end{array}$$
			\item Let $|u_2|>|v_2|$.  Then, there exists $\beta\in\N$, such that $$({\bf z_2}^{u_2}-{\bf z_2}^{v_2})^{\rm hom}={\bf z_2}^{u_2}-{\bf z_2}^{v_2}x^\beta.$$ Thus $|u_2|=|v_2|+\beta$. Furthermore, from the equation \ref{eq-aux1} and \ref{deg-hom}, we deduce that 
			$$(|u_1|+\alpha)+|u_2|=|v_1|+(|v_2|+\beta)=(|v_1|+|v_2|)+\beta=(|u_1|+|u_2|+\delta)+\beta$$
			So $\alpha=\delta+\beta$ and
			$$\begin{array}{ll}
			B&={\bf z_1}^{u_1}{\bf z_2}^{u_2}x^\delta-{\bf z_1}^{v_{1}}{\bf z_2}^{v_{2}}\\
			&={\bf z_1}^{u_1}x^\delta({\bf z_2}^{u_2}-{\bf z_2}^{v_2}x^\beta)+{\bf z_2}^{v_2}({\bf z_1}^{u_1}x^{\alpha}-{\bf z_1}^{v_1})\\
			&={\bf z_1}^{u_1}x^\delta({\bf z_2}^{u_2}-{\bf z_2}^{v_2})^{\rm hom}+{\bf z_2}^{v_2}({\bf z_1}^{u_1}-{\bf z_1}^{v_1})^{{\rm hom}}\in I_2+I_1.
			\end{array}$$
		\end{enumerate}
	 \item Let $|u_1|>|v_1|$. Then, from the equation \ref{deg-hom} $|v_2|>|u_2|$. By this manner, there exist $\alpha,\beta\in\N$, $\alpha>0$ and $\beta>0$ such that
		$$({\bf z_1}^{u_1}-{\bf z_1}^{v_1})^{\rm hom}={\bf z_1}^{u_1}-{\bf z_1}^{v_1}x^{\alpha}\mbox{ }{\rm et}\mbox{ }({\bf z_2}^{u_2}-{\bf z_2}^{v_2})^{\rm hom}={\bf z_2}^{u_2}x^\beta-{\bf z_2}^{v_2}.$$
	So $|u_1|+(|u_2|+\beta)=(|v_1|+\alpha)+|v_2|=(|v_1|+|v_2|)+\alpha=(|u_1|+|u_2|+\delta)+\alpha$, this last equality is due to the equation \ref{deg-hom}. By consequence $\beta=\delta+\alpha$ and finally we have:
	$$\begin{array}{ll}
			B&={\bf z_1}^{u_1}{\bf z_2}^{u_2}x^\delta-{\bf z_1}^{v_{1}}{\bf z_2}^{v_{2}}\\
			&={\bf z_2}^{u_2}x^\delta({\bf z_1}^{u_1}-{\bf z_1}^{v_1}x^{\alpha})+{\bf z_2}^{v_1}({\bf z_2}^{u_2}x^\beta-{\bf z_2}^{v_2})\\
			&={\bf z_2}^{u_2}x^\delta({\bf z_1}^{u_1}-{\bf z_1}^{v_1})^{{\rm hom}}+{\bf z_2}^{v_1}({\bf z_2}^{u_2}x^\beta-{\bf z_2}^{v_2})^{{\rm hom}}\in I_1+I_2.
	\end{array}$$
	 \item Let $|u_1|=|v_1|$. Then, from the equation \ref{deg-hom}, $|v_2|\geq|u_2|$. By this way, there exist $\beta\in\N$, such that
		$$({\bf z_1}^{u_1}-{\bf z_1}^{v_1})^{\rm hom}={\bf z_1}^{u_1}-{\bf z_1}^{v_1}\mbox{ }{\rm et}\mbox{ }({\bf z_2}^{u_2}-{\bf z_2}^{v_2})^{\rm hom}={\bf z_2}^{u_2}x^\beta-{\bf z_2}^{v_2}.$$
Thus $|u_1|+(|u_2|+\beta)=|v_1|+|v_2|=|v_1|+|v_2|=|u_1|+|u_2|+\delta$, this last equality is due to the equation \ref{deg-hom}. By consequence $\beta=\delta$ and we have:
	$$\begin{array}{ll}
			B&={\bf z_1}^{u_1}{\bf z_2}^{u_2}x^\delta-{\bf z_1}^{v_{1}}{\bf z_2}^{v_{2}}\\
			&={\bf z_2}^{u_2}x^\delta({\bf z_1}^{u_1}-{\bf z_1}^{v_1})+{\bf z_2}^{v_1}({\bf z_2}^{u_2}x^\beta-{\bf z_2}^{v_2})\\
			&={\bf z_2}^{u_2}x^\delta({\bf z_1}^{u_1}-{\bf z_1}^{v_1})^{{\rm hom}}+{\bf z_2}^{v_1}({\bf z_2}^{u_2}x^\beta-{\bf z_2}^{v_2})^{{\rm hom}}\in I_1+I_2.
	\end{array}$$
	\end{enumerate}
	\end{itemize}
It proves that $J\subset I_1+I_2$ and we conclude $I_1+I_2=J:={\rm Ker}\mbox{ }\overline{\varphi}$.
\end{itemize}
Finally, by the remark \ref{remarque-def-torique} we conclude $I_1+I_2$ is a toric homogeneous ideal.

\item We have proved that the matrix $A\in M_{m_1+m_2+1,n_1+n_2+1}(\Z)$ which represents the parametrization $\overline{\varphi}$ of $I_1+I_2$ is given by:
$$A=\left(\begin{array}{c|c|c}
A_1'&{\bf 0}&{\bf 0}\\ \hline
{\bf 0}&A_2'&{\bf 0}\\ \hline
\alpha_1&\alpha_2&\gamma
\end{array}\right).$$
By consequence: ${\rm dim}\mbox{ }(\k[{\bf z_{1}},{\bf z_{2}},x]/(I_1+I_2))={\rm dim}\mbox{ } (\k[{\bf z_{1}},x]/I_1)+{\rm dim}\mbox{ }(\k[{\bf z_{2}},x]/I_2)-1$.
\end{enumerate}

\qed

\igno{\begin{next}\nx{Exemple}\label{no-torique}
 
\end{next}
}

\begin{next}\nx{Proposition}\label{torique1}
Let ${\bf z_{1}},\dots,{\bf z_{k}}$ be $k$ disjoint sets of variables and for any $i=1,\dots, k$, let $I_i$ be a homogeneous toric ideal with the maximal rank parametrizations $\varphi_i:\k[{\bf z_i}]\rightarrow\k[{\bf t_i}^{\pm 1}]$. Then $\sum_{i=1}^{k}I_i$ is a homogeneous toric ideal with maximal rank parametrization
\begin{center}
$\varphi:\k[{\bf z_{1}},\dots,{\bf z_{k}}]\rightarrow \k[{\bf t_1}^{\pm 1},\dots,{\bf t_k}^{\pm 1}]$ the morphism defined by $\varphi(z_{i,j})=\varphi_i(z_{i,j})$.
\end{center} 
Moreover
$${\rm dim}\mbox{ }(k[{\bf z_{1}},\dots,{\bf z_{k}}]/(\sum_{i=1}^{k}I_i))=\sum_{i=1}^{k}{\rm dim}\mbox{ }(\k[{\bf z_{i}}]/I_i).$$
\end{next}
\proof
We set $J:=\sum_{i=1}^{k}I_i$. It is clear that $J\subset {\rm Ker}\mbox{ }\varphi$. We will prove that ${\rm Ker}\mbox{ }\varphi\subset J$:\\
Let $m_1,m_2$ be monomials of $\k[{\bf z_{1}},\dots,{\bf z_{k}}]$ such that $b=m_1-m_2$ is a generator of ${\rm Ker}\mbox{ }\varphi$. Furthermore, we write $m_1=m_1^{(1)}m_1^{(2)}\cdots m_1^{(k)}$ and $m_2=m_2^{(1)}m_2^{(2)}\cdots m_2^{(k)}$, where for any $1\leq i\leq k$ $m_1^{(i)}$ and $m_2^{(i)}$ are monomials in $\k[{\bf z_{i}}]$. Since $\varphi(m_1)=\varphi(m_2)$ and ${\bf t_1},\dots,{\bf t_k}$ are disjoint sets of variable,\break $\varphi(m_1^{(i)})=\varphi(m_2^{(i)})$. So, for any $1\leq i<j\leq k$:  $$\varphi_i(m_1^{(i)})=\varphi(m_1^{(i)})=\varphi(m_2^{(i)})=\varphi_i(m_2^{(i)}),$$ thus $m_1^{(i)}-m_2^{(i)}\in {\rm Ker}\mbox{ }\varphi_i=\I_{i}$ and
$$m_1-m_2=m_1^{(2)}\cdots m_1^{(k)}(m_1^{(1)}-m_2^{(1)})+m_2^{(1)}m_1^{(3)}\cdots m_1^{(k)}(m_1^{(2)}-m_2^{(2)})\igno{+m_2^{(1)}m_2^{(2)}m_1^{(4)}\cdots m_1^{(k)}(m_1^{(3)}-m_2^{(3)})}+\dots +m_2^{(1)}\cdots m_2^{(k-1)}(m_1^{(k)}-m_2^{(k)})\in J.$$
Thus ${\rm Ker}\mbox{ }\varphi\subset J$ and we conclude that $J={\rm Ker}\mbox{ }\varphi$. So $J$ is a toric ideal. It is clear that the parametrizatrion $\varphi$ of $J$ is a maximal rang parametrization, thus $${\rm dim}\mbox{ }(k[{\bf z_{1}},\dots,{\bf z_{k}}]/(\sum_{i=1}^{k}I_i))=\sum_{i=1}^{k}{\rm dim}\mbox{ }(\k[{\bf z_{i}}]/I_i).$$
\qed

\begin{next}\nx{Definition}\label{some-scroll-torique}
{\rm Let ${\bf x_1},\dots,{\bf x_k}$ be variable sets and $I_1\subset\k[{\bf x_1}],\dots,I_k\subset\k[{\bf x_k}]$ be toric ideals such that for all pair $i,j\in\{1,\dots,k\}$, $i\neq j$, $|{\bf x_i}\cap{\bf x_j}|\leq 1$. We define the {\bf graph of the sequence of toric ideals} $I_1,\dots,I_k$, denoted by $G(I_1,\dots,I_k)$,as the graph whose vertice set is $I_1,\dots,I_k$ and $\{I_i,I_j\}$ is an edge if $I_i$ and $I_j$ has a common variable.}
\end{next}

\begin{next}\nx{Theorem}\label{toriquin}
Let ${\bf x_1},\dots,{\bf x_k}$ be set of variables and $I_1\subset\k[{\bf x_1}],\dots,I_k\subset\k[{\bf x_k}]$ be toric ideals such that for any pair $i,j\in\{1,\dots,k\}$, $i\neq j$, $|{\bf x_i}\cap{\bf x_j}|\leq 1$. If for any connected component $G(I_1,\dots,I_k)$ is a tree, then $J:=\sum_{i=1}^k I_i$ is a homogeneous toric ideal and setting $r$ as the number of connected components of $G(I_1,\dots,I_k)$ we have that:
$${\rm dim}\mbox{ }(\si/J)=\sum_{i=1}^k{\rm dim}\mbox{ }(\k[{\bf x_i}]/I_i)+r-k+1,$$
where ${\bf x}=\cup_{i=1}^k{\bf x_i}$.
\end{next}
\proof
Let $G_j$ be a connected component of $G:=G(I_1,\dots,I_k)$ and $V(G_j)=\{I_{j_1},\dots,I_{j_{s_j}}\}$ be the vertex set of $G_j$. We remark that $J:=\sum_{i=1}^k I_i=\sum_{j=1}^r(\sum_{i=1}^{s_j} I_{j_i})$. So, if we prove that for any $j\in\{1,\dots, r\}$, $(\sum_{i=1}^{s_j} I_{j_i})$ is a homogeneous toric ideal and $${\rm dim}\mbox{ }(\si/\sum_{i=1}^{s_j} I_{j_i})=\sum_{i=1}^{s_j}{\rm dim}\mbox{ }(\k[{\bf x_{j_i}}]/I_i)-s_j,$$ then $J$ is a homogeneous toric ideal due to the proposition \ref{torique1} and $${\rm dim}\mbox{ }(\si/J)=\sum_{i=1}^k{\rm dim}\mbox{ }(\k[{\bf x_i}]/I_i)+r-k+1.$$ Then, we can assume that $G$ is a connected tree and we will prove by induction over $k$ that $J$ is a homogenous toric ideal.
\begin{enumerate}
 \item If $k=2$, it is the theorem \ref{sum-2-torique}
 \item We assume that if $G$ is a connected tree with $k-1$-vertices, then $\sum_{i=1}^{k-1} I_i$ is a toric homogeneous ideal and $${\rm dim}\mbox{ }(\si/\sum_{i=1}^{k-1} I_{i})=\sum_{i=1}^{k-1}{\rm dim}\mbox{ }(\k[{\bf x_{i}}]/I_i)-(k-1).$$ We must prove that if $G$ is a connected tree with $k$-vertice, then $J:=\sum_{i=1}^k I_i$ is a toric homogeneous ideal and $${\rm dim}\mbox{ }(\si/J)=\sum_{i=1}^{k}{\rm dim}\mbox{ }(\k[{\bf x_{i}}]/I_i)-k.$$
It will follow from the next argument: since $G$ is a connected tree, there exists $i\in\{1,\dots,k\}$ such that there exists a unique $j\in\{1,\dots,k\}$ such that $\{I_i,I_j\}\in E(G)$. We can assume that $i=k$, thus $|{\bf x_i}\cap{\bf x_k}|=1$ and for any $j\in(\{1,\dots,k-1\}\setminus\{i\})$, $|{\bf x_{j}}\cap{\bf x_k}|=0$. So $|(\cup_{i=1}^{k-1}{\bf x_i})\cap{\bf x_k}|=1$ and $G\setminus\{I_k\}$ is a connected tree. By induction hypothesis $\sum_{i=1}^{k-1} I_i$ is a homogeneous toric ideal and $${\rm dim}\mbox{ }(\si/\sum_{i=1}^{k-1} I_{i})=\sum_{i=1}^{k-1}{\rm dim}\mbox{ }(\k[{\bf x_{i}}]/I_i)-(k-1).$$ In addition, $I_k$ is a homogeneous toric ideal, so by the theorem \ref{sum-2-torique} $J=\sum_{i=1}^{k-1} I_i+I_k$ is a homogeneous toric ideal and 
$$\begin{array}{ll}
   {\rm dim}\mbox{ }(\si/J)&={\rm dim}\mbox{ }(\si/\sum_{i=1}^{k-1} I_{i})+{\rm dim}\mbox{ }(\k[{\bf x_k}]/I_{k})\\
	&=\sum_{i=1}^{k-1}{\rm dim}\mbox{ }(\k[{\bf x_{i}}]/I_i)-(k-1)+ {\rm dim}\mbox{ }(\k[{\bf x_{k}}]/I_k)-1\\
	&=\sum_{i=1}^{k}{\rm dim}\mbox{ }(\k[{\bf x_{i}}]/I_i)-k.
  \end{array}
$$

\end{enumerate}

\qed

\end{document}